# ON THE MAXIMUM BIAS FUNCTIONS OF MM-ESTIMATES AND CONSTRAINED M-ESTIMATES OF REGRESSION

BY JOSÉ R. BERRENDERO[1], BEATRIZ V. M. MENDES[2] AND
DAVID E. TYLER[3]

*Universidad Autónoma de Madrid, Federal University of Rio de Janeiro
and Rutgers University*

We derive the maximum bias functions of the *MM*-estimates and the constrained *M*-estimates or *CM*-estimates of regression and compare them to the maximum bias functions of the *S*-estimates and the $\tau$-estimates of regression. In these comparisons, the *CM*-estimates tend to exhibit the most favorable bias-robustness properties. Also, under the Gaussian model, it is shown how one can construct a *CM*-estimate which has a smaller maximum bias function than a given *S*-estimate, that is, the resulting *CM*-estimate *dominates* the *S*-estimate in terms of maxbias and, at the same time, is considerably more efficient.

**1. Introduction.** An important consideration for any estimate is an understanding of its robustness properties. Different measures exist which try to reflect the general concept known as robustness. One such measure is the maximum bias function, which measures the maximum possible bias of an estimate under $\varepsilon$-contamination. In this paper, we study the maximum bias functions for the *MM*-estimates and the constrained *M*-estimates or *CM*-estimates of regression and compare them to the maximum bias functions for the *S*-estimates and the $\tau$-estimates of regression.

The maximum bias functions for Rousseeuw and Yohai's [10] *S*-estimates of regression were originally derived by Martin, Yohai and Zamar [7] under the assumption that the independent variables follow an elliptical distribution and that the intercept term is known. More recently, Berrendero

Received February 2006.
[1]Supported in part by Spanish Grant MTM2004-00098 and Grant S-0505/ESP/0158 of the Comunidad de Madrid.
[2]Supported in part by CNPq-Brazil, Edital CNPq 019/2004 Universal.
[3]Supported in part by NSF Grant DMS-03-05858.

and Zamar [1] derived the maximum bias functions for the $S$-estimates of regression under much broader conditions. Further general results on the maximum bias functions can be found in [4]. The method used in [1] applies to a wide class of regression estimates. For example, it allows one to obtain the maximum bias functions of Yohai and Zamar's [12] $\tau$-estimates of regression. Unfortunately, it does not apply to Yohai's [11] *MM*-estimates of regression, arguably the most popular high breakdown point estimates of regression. The *MM*-estimates, for example, are the default robust regression estimates in S-PLUS.

The original motivation for the current paper was thus to derive the maximum bias functions of the *MM*-estimates of regression and compare them to the maximum bias functions of the $S$-estimates and $\tau$-estimates of regression. A lesser known high breakdown point estimate of regression, namely Mendes and Tyler's [8] constrained $M$-estimates of regression (or *CM*-estimates for short), has also been included in the study since the associated maximum bias functions can be readily obtained by applying the general method given in [1]. Expressions for the maximum bias functions of the *MM*-estimates and the *CM*-estimates are derived in Sections 3 and 4. Comparisons between the $S$-, $\tau$-, *MM*- and *CM*-estimates based on biweight score functions are given in Section 5. It turns out that in these comparisons, the *CM*-estimates tend to exhibit the most favorable robustness properties.

Consequently, a more detailed theoretical comparison between the maximum bias functions of the $S$-estimates and the *CM*-estimates of regression, which helps explain the computational comparisons made in Section 5, is given in Section 6. In particular, under the Gaussian model, it is shown how one can construct a *CM*-estimate of regression so that its maximum bias function *dominates* that of a given $S$-estimate of regression. That is, the maximum bias function of the *CM*-estimate is smaller for some level of contamination $\varepsilon$ and is never larger for any value of $\varepsilon$. The $S$-estimate is thus said to be *bias-inadmissible* at the Gaussian model.

Section 2 reviews the notion of the maximum bias function in the regression setting, as well as the definitions of the $S$-estimates, the *MM*-estimates and the *CM*-estimates for regression. Technical proofs are given in the Appendix.

**2. The regression model and the concept of maximum bias.** We follow the general setup given in [7]. Specifically, we consider the linear regression model

$$(2.1) \qquad y = \alpha_o + \mathbf{x}'\boldsymbol{\theta}_o + u,$$

where $y \in \mathbb{R}$ represents the response, $\mathbf{x} = (x_1, x_2, \ldots, x_p)' \in \mathbb{R}^p$ represents a random vector of explanatory variables, $\alpha_o \in \mathbb{R}$ and $\boldsymbol{\theta}_o \in \mathbb{R}^p$ are the true intercept and slope parameters, respectively, and the random error term $u \in \mathbb{R}$ is assumed to be independent of $\mathbf{x}$. Let $F_o$ and $G_o$ represent the distribution functions of $u$ and $\mathbf{x}$, respectively, and let $H_o$ represent the corresponding joint distribution function of $(y, \mathbf{x})$. The following assumptions on the distribution $H_o$ are assumed throughout the paper:

A1. $F_o$ is absolutely continuous with density $f_o$ which is symmetric, continuous and strictly decreasing on $\mathbb{R}^+$;

A2. $P_{G_o}(\mathbf{x}'\boldsymbol{\theta} = c) < 1$ for any $\boldsymbol{\theta} \in \mathbb{R}^p$, $\boldsymbol{\theta} \neq \mathbf{0}$, $c \in \mathbb{R}$.

As in [7] and [1], we focus on the estimation of the slope parameters $\boldsymbol{\theta}_o$. One reason for doing so is that once given a good estimate of the slope parameters, the problem of estimating the intercept term and the residual scale reduces to the well-studied univariate location and scale problem. Let $\mathbf{T}$ represent some $\mathbb{R}^p$-valued functional defined on $\mathcal{H}$, a space of distribution functions on $\mathbb{R}^{p+1}$ which includes some weak neighborhood of $H_o$, such that $\mathbf{T}(H_o) = \boldsymbol{\theta}_o$. For sufficiently large $n$, $\mathcal{H}$ almost surely contains the empirical distribution function $H_n$ corresponding to a random sample $\{(y_1, \mathbf{x}_1), \ldots, (y_n, \mathbf{x}_n)\}$ from $H_o$. Furthermore, we assume that $\mathbf{T}$ is weakly continuous at $H_o$, and so the statistic $\mathbf{T}_n = \mathbf{T}(H_n)$ is a consistent estimate of $\boldsymbol{\theta}_o$.

All functionals $\mathbf{T}$ considered in this paper are regression equivariant, as defined, for example, in [7]. For such functionals, a natural invariant measure of the "asymptotic" bias of $\mathbf{T}$ at $H$ is given by

$$(2.2) \quad b_{\boldsymbol{\Sigma}_o}(\mathbf{T}, H) = \begin{cases} ((\mathbf{T}(H) - \boldsymbol{\theta}_o)' \boldsymbol{\Sigma}_o (\mathbf{T}(H) - \boldsymbol{\theta}_o))^{1/2}, & H \in \mathcal{H}, \\ \infty, & H \notin \mathcal{H}. \end{cases}$$

Here, $\boldsymbol{\Sigma}_o = \boldsymbol{\Sigma}(G_o)$ is taken to be an affine equivariant scatter matrix for the regressors $\mathbf{x}$ under $G_o$. We can thus presume, without loss of generality, that $(\alpha_o, \boldsymbol{\theta}_o) = \mathbf{0}$ and $\boldsymbol{\Sigma}_o = \mathbf{I}$. Hence, the asymptotic bias of $\mathbf{T}$ at $H$ becomes the Euclidean norm of $\mathbf{T}$,

$$(2.3) \quad b(\mathbf{T}, H) = \begin{cases} \|\mathbf{T}(H)\|, & H \in \mathcal{H}, \\ \infty, & H \notin \mathcal{H}, \end{cases}$$

where $\mathcal{H}$ is the class of distributions such that $\|\mathbf{T}(H)\| < \infty$. The maximum asymptotic bias of $\mathbf{T}$ over $\varepsilon$-contaminated neighborhoods $V_\varepsilon$ of $H_o$, that is, $V_\varepsilon = \{H \mid H = (1 - \varepsilon)H_o + \varepsilon H^*, H^* \in \mathcal{H}^*\}$, where $\mathcal{H}^*$ is the set of all distribution functions on $\mathbb{R}^{p+1}$, is defined to be

$$(2.4) \quad B_{\mathbf{T}}(\varepsilon) = \sup\{b(\mathbf{T}, H) \mid H \in V_\varepsilon\}$$

and the asymptotic breakdown point is subsequently defined to be

$$(2.5) \quad \varepsilon^* = \inf\{\varepsilon \mid B_{\mathbf{T}}(\varepsilon) = \infty\}.$$

From an applied perspective, regardless of $\boldsymbol{\Sigma}_o$, it may be of interest to derive upper bounds for the Euclidean distance between $\mathbf{T}(H)$ and $\boldsymbol{\theta}_o$, that is, for $\|\mathbf{T}(H)-\boldsymbol{\theta}_o\|$. This measure is referred to as a *bias bound* by Berrendero and Zamar in [1], wherein they use it for adjusting confidence intervals for $\boldsymbol{\theta}$ to include the possibility of bias introduced by a contaminated model. Note that the *bias bound* is regression and scale equivariant, but not affine equivariant, and hence is not directly related to the maximum bias (2.4). In [1] some results are given for computing *bias bounds*, taking the maximum bias function as a starting point.

2.1. *M-estimates with general scale.* The *S*-, *MM*- and *CM*-estimates of regression all lie within the class of $M$-estimates with general scale considered in [7]. An $M$-estimate, or, more appropriately, an $M$-functional, with general scale for the regression parameters $\alpha_o$ and $\boldsymbol{\theta}_o$, say $\mathrm{t}(H)$ and $\mathbf{T}(H)$, respectively, can be defined as the solution which minimizes

$$\mathrm{E}_H\left[\rho\left(\frac{y-\alpha-\mathbf{x}'\boldsymbol{\theta}}{\sigma(H)}\right)\right] \tag{2.6}$$

over all $\alpha \in \mathbb{R}$ and $\boldsymbol{\theta} \in \mathbb{R}^p$, where $\rho$ is some nonnegative symmetric function and $\sigma(H)$ is some scale functional. The scale functional $\sigma(H)$ may be determined simultaneously or independently of $\{\mathrm{t}(H),\mathbf{T}(H)\}$. We assume throughout the paper that $\sigma(H)$ is regression invariant and residual scale equivariant, again as defined, for example, in [7]. Throughout, it is assumed that the function $\rho$ satisfies the following conditions:

A3.  (i) $\rho$ is symmetric and nondecreasing on $[0,\infty)$ with $\rho(0)=0$;
   (ii) $\rho$ is bounded with $\lim_{u\to\infty}\rho(u)=1$;
   (iii) $\rho$ has only a finite number of discontinuities.

If the function $\rho$ is also differentiable, then $\{\mathrm{t}(H),\mathbf{T}(H)\}$ is a solution to the $p+1$ simultaneous $M$-estimating equations

$$(2.7)\ \mathrm{E}_H\left\{\psi\left(\frac{y-\alpha-\mathbf{x}'\boldsymbol{\theta}}{\sigma(H)}\right)\mathbf{x}\right\}=\mathbf{0} \quad \text{and} \quad \mathrm{E}_H\left\{\psi\left(\frac{y-\alpha-\mathbf{x}'\boldsymbol{\theta}}{\sigma(H)}\right)\right\}=0,$$

where $\psi(u) \propto \rho'(u)$. By Condition A3(i), $\psi$ is an odd function, nonnegative on $[0,\infty)$. Condition A3(ii) implies that these $M$-estimates are redescending, that is, $\psi(u) \to 0$ as $u \to \infty$. A popular choice for $M$-estimates are Tukey's biweighted $M$-estimates, which correspond to choosing $\rho(u)$ to be

$$\rho_T(u) = \begin{cases} 3u^2 - 3u^4 + u^6, & \text{for } |u| \le 1, \\ 1, & \text{for } |u| > 1. \end{cases} \tag{2.8}$$

Note that this gives rise to the biweight $\psi$ function $\psi_T(u) = u\{(1-u^2)_+\}^2$.

The $S$-estimates for the intercept, slopes and scale are collectively defined to be the solution $\{t_s(H), \mathbf{T}_s(H), \sigma_s(H)\}$ to the problem of minimizing $\sigma \in \mathbb{R}^+$ subject to the constraint

$$\text{(2.9)} \qquad \text{E}_H\left[\rho\left(\frac{y - \alpha - \mathbf{x}'\boldsymbol{\theta}}{\sigma}\right)\right] \leq b$$

for some fixed value $b$, $0 < b < 1$. The breakdown point of the $S$-estimate of regression is $\varepsilon^* = \min\{b, 1-b\}$. A drawback to the $S$-estimates is that the tuning constant $b$ not only determines the breakdown point, but it also determines the efficiency of the estimate. To obtain a reasonable efficiency under a normal error model, one must usually decrease the breakdown point substantially.

This problem with tuning the $S$-estimates of regression motivated Yohai [11] to introduce the *MM*-estimates of regression which can be tuned to have high efficiency under normal error while simultaneously maintaining a high breakdown point. Let $\rho_1$ and $\rho_2$ be a pair of loss functions satisfying A3 and with $\rho_1 > \rho_2$. Set $b = \text{E}_{F_o}\rho_1(Y)$. *MM*-estimates are collectively defined to be the solution $\{t_{MM}(H), \mathbf{T}_{MM}(H)\}$ which minimizes

$$L_H(\alpha, \boldsymbol{\theta}) = \text{E}_H\left[\rho_2\left(\frac{y - \alpha - \mathbf{x}'\boldsymbol{\theta}}{s(H)}\right)\right],$$

where $s(H) \doteq \sigma_s(H)$ is a preliminary $S$-functional of scale, as defined above, based on $\rho = \rho_1$. The breakdown point of the *MM*-estimates depends only on $\rho_1$ and is given by $\varepsilon^* = \min\{b, 1-b\}$. On the other hand, their asymptotic distribution is determined exclusively by $\rho_2$. This allows the *MM*-estimates to be tuned so that they possess both high breakdown point and high efficiency.

The *CM*-estimates are another class of regression estimates which can be tuned to have high efficiency at the normal model while maintaining a high breakdown point. The *CM*-estimates for the intercept, slopes and scale are collectively defined to be the solution $\{t_{CM}(H), \mathbf{T}_{CM}(H), \sigma_{CM}(H)\}$ which minimizes

$$\text{(2.10)} \qquad \text{L}_H(\alpha, \boldsymbol{\theta}, \sigma) = c\text{E}_H\left[\rho\left(\frac{y - \alpha - \mathbf{x}'\boldsymbol{\theta}}{\sigma}\right)\right] + \log \sigma$$

subject to the constraint (2.9), where $c > 0$ represents a tuning constant. As with the $S$-estimates of regression, the asymptotic breakdown point of the *CM*-estimates of regression is $\varepsilon^* = \min\{b, 1-b\}$. Unlike the $S$-estimates of regression, though, the *CM*-estimates of regression can be tuned by means of the constant $c$ in order to obtain a reasonably high efficiency without affecting the breakdown point.

We again emphasize that our focus here is on the slope functionals $\mathbf{T}(H)$, rather than on the intercept functionals $t(H)$ or the scale functionals $\sigma(H)$.

Given a good slope functional, one may wish to consider the wider range of location and scale functionals based on the distribution of $y - \mathbf{x}'\mathbf{T}(H)$ such as its median and median absolute deviation, rather than those arising from an *S*-, *MM*- or *CM*-estimate of regression.

### 3. Maximum bias functions.

3.1. *Maximum bias functions for MM-estimates.* If $F_{H,\alpha,\boldsymbol{\theta}}$ is the distribution function of the absolute residuals $|y - \alpha - \mathbf{x}'\boldsymbol{\theta}|$, then Berrendero and Zamar [1] give an expression for the maximum bias function for any estimate whose definition can be expressed in the form

$$(3.11) \qquad \{t(H), \mathbf{T}(H)\} = \arg\min_{(\alpha,\boldsymbol{\theta})} J(F_{H,\alpha,\boldsymbol{\theta}}),$$

where $J(F)$ is a functional possessing certain monotonic properties. The *S*-, $\tau$- and *CM*-estimates are of this form. Applications of their general results to the *S*- and $\tau$-estimates are given in [1]. Application of these results to the *CM*-estimates is presented in Section 3.2.

The *MM*-estimates, however, cannot be expressed in the form (3.11), so a different approach is needed in order to study their bias behavior. Let $B_{MM}(\varepsilon)$ be the maximum bias function of an *MM*-estimate of regression. In this subsection, lower and upper bounds for $B_{MM}(\varepsilon)$ are obtained under quite general conditions. In some important cases, these two bounds are often equal and thereby allow for the exact determination of the maximum bias function.

Let $\underline{s} = \inf_{H \in V_\varepsilon} s(H)$, $\overline{s} = \sup_{H \in V_\varepsilon} s(H)$ and

$$m(t,s) = \inf_{\|\boldsymbol{\theta}\|=t} \inf_{\alpha \in \mathbb{R}} E_{H_o} \rho_2\left(\frac{y - \alpha - \mathbf{x}'\boldsymbol{\theta}}{s}\right) - E_{H_o} \rho_2\left(\frac{y}{s}\right).$$

The following two functions play a key role in the developments below:

$$(3.12) \qquad h_1(t) = m(t, \overline{s}) \quad \text{and} \quad h_2(t) = \inf_{\underline{s} \leq s \leq \overline{s}} m(t, s).$$

THEOREM 3.1. *Let $\mathbf{T}_{MM}$ be an MM-estimate of the regression slopes with loss functions $\rho_i$, $i = 1, 2$, satisfying* A3. *Assume that the maximum bias function of the S-estimate with score function $\rho_1$, $B_S(\varepsilon)$, satisfies $B_S(\varepsilon) < h_1^{-1}[\varepsilon/(1-\varepsilon)]$. Under* A1 *and* A2, *the maximum bias function of $\mathbf{T}_{MM}$, $B_{MM}(\varepsilon)$, satisfies*

$$(3.13) \qquad h_1^{-1}\left(\frac{\varepsilon}{1-\varepsilon}\right) \leq B_{MM}(\varepsilon) \leq h_2^{-1}\left(\frac{\varepsilon}{1-\varepsilon}\right).$$

Note that the condition $B_S(\varepsilon) < h_1^{-1}[\varepsilon/(1-\varepsilon)]$ of the above theorem, together with (3.13), implies that $B_S(\varepsilon) < B_{MM}(\varepsilon)$. This condition usually holds for an appropriately chosen $\rho_1$ function. Thus, an *MM*-estimate does not improve upon the maximum bias of the initial *S*-estimate. The trade-off, though, is that with an appropriately chosen $\rho_2$ function, the *MM*-estimate can greatly improve upon the efficiency of the initial *S*-estimate.

Upper and lower bounds for the maximum bias of *MM*-estimates have also been obtained, respectively, by Hennig [5], Theorem 3.1, and Martin, Yohai and Zamar [7], Lemma 4.1, under the assumption of unimodal elliptically distributed regressors. For this special case, the upper bounds given in (3.13) and in [5] agree. On the other hand, the lower bound given in [7], namely $B_{MM}(\varepsilon) \geq h_0^{-1}[\varepsilon/(1-\varepsilon)]$, where $h_0(t) = \sup_{\underline{s} \leq s \leq \overline{s}} m(t,s)$, is not as tight as that given in (3.13).

In our setup, the assumption of unimodal elliptical regressors is equivalent to the following:

A2*. Under $G_o$, the distribution of $\mathbf{x}'\boldsymbol{\theta}$ is absolutely continuous, with a symmetric, unimodal density, and depends on $\boldsymbol{\theta}$ only through $\|\boldsymbol{\theta}\|$ for all $\boldsymbol{\theta} \neq \mathbf{0}$.

Under this condition, we can define

$$(3.14) \qquad g(s,t) = \mathrm{E}_{H_o}\left[\rho\left(\frac{y - \mathbf{x}'\boldsymbol{\theta}}{s}\right)\right],$$

where $\boldsymbol{\theta}$ is any vector such that $\|\boldsymbol{\theta}\| = t$. Under conditions A1, A2* and A3, it is shown in Lemma 3.1 of Martin, Yohai and Zamar [7] that $g$ is continuous, strictly increasing with respect to $\|\boldsymbol{\theta}\|$ and strictly decreasing in $s$ for $s > 0$.

If A2* holds, then $\underline{s}$ and $\overline{s}$ are defined so that $g_1(\underline{s}, 0) = b/(1-\varepsilon)$ and $g_1(\overline{s}, 0) = (b-\varepsilon)/(1-\varepsilon)$, respectively, and $m(t,s) = g_2(s,t) - g_2(s,0)$, where $g_i(s,t)$ is defined as in (3.14) after replacing $\rho$ with $\rho_i$.

3.2. *Maximum bias curves for CM-estimates.* A *CM*-estimate of regression $\{\mathrm{t}_{CM}(H), \mathbf{T}_{CM}(H)\}$ can be expressed in the form (3.11) with $J$ taken to be

$$(3.15) \qquad J_{CM}(F) = \inf_{s \geq \sigma(F)} c\mathrm{E}_F[\rho(\mathrm{y}/s)] + \log s$$

and where $\sigma(F)$ is the *M*-scale defined as the solution to the equation

$$(3.16) \qquad \mathrm{E}_F[\rho(\mathrm{y}/\sigma(F))] = b.$$

Consequently, application of the general method in [1] for computing maximum bias functions leads to the following result:

THEOREM 3.2. *Let $\mathbf{T}_{CM}$ be a CM-estimate of the regression slopes based on a function $\rho$ satisfying* A3, *and suppose $H_o$ satisfies* A1 *and* A2. *Define*

$$r_{CM}(\varepsilon) = J_{CM}[(1-\varepsilon)F_{H_0,0,\mathbf{0}} + \varepsilon\delta_\infty]$$

*and let*

(3.17)  $$m_{CM}(t) = \inf_{\|\boldsymbol{\theta}\|=t} \inf_{\alpha\in\mathbb{R}} J_{CM}[(1-\varepsilon)F_{H_0,\alpha,\boldsymbol{\theta}} + \varepsilon\delta_0].$$

*Then the maximum bias function of $\mathbf{T}_{CM}$, denoted by $B_{CM}(\varepsilon)$, is given by*

(3.18)  $$B_{CM}(\varepsilon) = m_{CM}^{-1}[r_{CM}(\varepsilon)].$$

This general result can be given a simpler representation when condition A2* also holds. In particular, in the definition of $m_{CM}(t)$, the infimum is obtained when $\alpha = 0$ and $\boldsymbol{\theta}$ is any vector such that $\|\boldsymbol{\theta}\| = t$. This gives

$$m_{CM}(t) = \inf\{A_{c,\varepsilon}(s,t) \mid s \geq m_s(t)\},$$

where $A_{c,\varepsilon}(s,t) = c(1-\varepsilon)g(s,t) + \log s$ and $m_s(t) = g_{(1)}^{-1}(b/(1-\varepsilon),t)$, with $g(s,t)$ being defined as in (3.14) and $g_{(1)}^{-1}(\cdot,t)$ being the inverse of $g$ with respect to $s$. Also, it is easy to verify that

$$r_{CM}(\varepsilon) = \inf\{A_{c,\varepsilon}(s,0) \mid s \geq r_s(\varepsilon)\} + c\varepsilon,$$

where $r_s(\varepsilon) = g_{(1)}^{-1}((b-\varepsilon)/(1-\varepsilon),0)$.

**4. Maximum bias functions for two special cases.** Maximum bias functions generally tend to have rather complicated expressions. Under some model distributions, though, these expressions can be substantially simplified. This is possible for two special cases considered here, namely the Gaussian and Cauchy models. These simplified expressions are useful for computing and comparing the maximum bias curves of various estimates for these models, which is done in Section 5.

4.1. *Maximum bias functions under the Gaussian model.* We assume throughout this section that not only the error term, but also that the regressor variables arise from a multivariate normal distribution. That is, we assume $H_o$ has a joint $N(\mathbf{0}, \mathbf{I}_{p+1})$ distribution and refer to this as the Gaussian model. Let $g(s) = E_\Phi \rho(Z/s)$, where $Z$ is a standard normal random variable, and define $\sigma_{b,\varepsilon} \doteq g^{-1}[(b-\varepsilon)/(1-\varepsilon)]$ and $\gamma_{b,\varepsilon} \doteq g^{-1}[b/(1-\varepsilon)]$. Martin, Yohai and Zamar [7] show that the maximum bias function for an $S$-estimate of the regression slope under the Gaussian model and based on a function $\rho$ satisfying A3 is given by

(4.19)  $$B_S(\varepsilon) = \left[\left(\frac{\sigma_{b,\varepsilon}}{\gamma_{b,\varepsilon}}\right)^2 - 1\right]^{1/2}.$$

To obtain an expression for the maximum bias function of a *CM*-estimate of regression under the Gaussian model, let

$$A_{c,\varepsilon}(s) = c(1-\varepsilon)g(s) + \log s. \tag{4.20}$$

Also, define $D_c(\varepsilon) = \inf_{s \geq \sigma_{b,\varepsilon}} A_{c,\varepsilon}(s) - \inf_{s \geq \gamma_{b,\varepsilon}} A_{c,\varepsilon}(s)$. We then have the following result:

THEOREM 4.1. *Let* $\mathbf{T}_{CM}$ *be a CM-estimate of the regression slopes based on a function $\rho$ satisfying* A3 *and assume* $H_o$ *is multivariate normal. We then have*

$$B_{CM}(\varepsilon) = \{\exp[2c\varepsilon + 2D_c(\varepsilon)] - 1\}^{1/2}. \tag{4.21}$$

Turning now to the *MM*-estimates, let $g_i(s) = E_\Phi \rho_i(Z/s)$ for $i = 1, 2$, where $Z$ is a standard normal random variable. Under the Gaussian model,

$$m(t,s) = g_2\left(\frac{s}{(1+t^2)^{1/2}}\right) - g_2(s).$$

Moreover, $\overline{s} = g_1^{-1}[(b-\varepsilon)/(1-\varepsilon)]$ and $\underline{s} = g_1^{-1}[b/(1-\varepsilon)]$. Since $\rho_1$ is the same $\rho$-function used in defining the preliminary *S*-estimate, we have $\overline{s} = \sigma_{b,\varepsilon}$ and $\underline{s} = \gamma_{b,\varepsilon}$. Hence, $B_{MM}(\varepsilon) \geq \ell(\varepsilon)$, where

$$\ell(\varepsilon) = h_1^{-1}\left(\frac{\varepsilon}{1-\varepsilon}\right)$$
$$= \left[\left(\frac{\sigma_{b,\varepsilon}}{g_2^{-1}[g_2(\sigma_{b,\varepsilon}) + \varepsilon/(1-\varepsilon)]}\right)^2 - 1\right]^{1/2}. \tag{4.22}$$

A simpler form for the upper bound, which can be used for computational purposes, can be obtained under some additional regularity conditions on $g_2(t)$. These conditions hold in most cases of interest.

A4. (i) $g(s)$ is continuously differentiable;
  (ii) $\phi(s) \doteq -sg'(s)$ is unimodal, with its maximum being obtained at $\sigma_M$. Set $K \doteq \phi(\sigma_M)$.

THEOREM 4.2. *In addition to the assumptions of Theorem* 3.1, *suppose that $g_2(s)$ satisfies* A4. *Then when $H_o$ is multivariate normal,*

$$\ell(\varepsilon) \leq B_{MM}(\varepsilon) \leq \max\{\ell(\varepsilon), u(\varepsilon)\},$$

*where $\ell(\varepsilon)$ is given in* (4.22) *and*

$$u(\varepsilon) = \left[\left(\frac{\gamma_{b,\varepsilon}}{g_2^{-1}[g_2(\gamma_{b,\varepsilon}) + \varepsilon/(1-\varepsilon)]}\right)^2 - 1\right]^{1/2}.$$

The upper bound in Theorem 4.2 coincides with that obtained by Hennig [5]. However, the tighter lower bound gives us further insight into the maximum bias and enables us to determine when the bounds are actually an equality. Obviously, if $\varepsilon$ is such that $u(\varepsilon) \leq \ell(\varepsilon)$, then $B_{MM}(\varepsilon) = \ell(\varepsilon)$. This occurs in many important cases for a wide range of $\varepsilon$ values.

As an example, consider the biweight loss function $\rho_T$ defined by (2.8). If we choose $\rho_1(u) = \rho_T(u/k_1)$ and $\rho_2(u) = \rho_T(u/k_2)$ with tuning constants $k_1 = 1.56$ and $k_2 = 4.68$, respectively, and choose $b = 0.5$, then the resulting $MM$-estimate has a 50% breakdown point and is asymptotically 95% efficient under the Gaussian model. For this case, it can be verified that the condition $B_S(\varepsilon) < h_1^{-1}[\varepsilon/(1-\varepsilon)]$ in Theorem 3.13 holds. From (4.22), it can be noted that this condition is equivalent to $g_2(\gamma_{b,\varepsilon}) - g_2(\sigma_{b,\varepsilon}) < \varepsilon/(1-\varepsilon)$. It can also be verified that the corresponding $\phi_2$ function is unimodal. A plot of $\phi_2$ is displayed in the upper panel of Figure 1. The bounds given in Theorem 4.2 for this $MM$-estimate are displayed in the lower panel of Figure 1. Both bounds coincide for values of $\varepsilon$ up to approximately 0.33. Hence, for such $\varepsilon$, the exact value at the maximum bias function is known.

4.2. *Maximum bias functions under the Cauchy model.* We now assume that the error term and the regressors follow independent Cauchy distributions rather than normal distributions. That is, we assume $\mathbf{x}_1, \ldots, \mathbf{x}_n$ and $y$ have independent standard Cauchy distributions. For brevity, we refer to this distributional model as the Cauchy model. Note that in this case, the distribution of the regressors is not elliptically symmetric. The derivations for the Cauchy model closely follow those given for the Gaussian model.

Let $g(s) = \mathrm{E}_\Phi \rho(Z/s)$, where $Z$ is now a standard Cauchy random variable, and, again, let $\sigma_{b,\varepsilon} \doteq g^{-1}[(b-\varepsilon)/(1-\varepsilon)]$ and $\gamma_{b,\varepsilon} \doteq g^{-1}[b/(1-\varepsilon)]$. In the Appendix, we show the maximum bias function for an $S$-estimate of regression to be

$$B_S(\varepsilon) = \frac{\sigma_{b,\varepsilon}}{\gamma_{b,\varepsilon}} - 1 \tag{4.23}$$

and for a $CM$-estimate of regression to be

$$B_{CM}(\varepsilon) = \exp\{D_c(\varepsilon) + c\varepsilon\} - 1, \tag{4.24}$$

with $D_c(\varepsilon)$ being analogous to its definition given after equation (4.20). Upper and lower bounds for the maximum bias function for the $MM$-estimates of regression are shown in the Appendix to be

$$\ell(\varepsilon) \leq B_{MM}(\varepsilon) \leq \max\{\ell(\varepsilon), u(\varepsilon)\} \quad \text{where} \tag{4.25}$$

$$\ell(\varepsilon) = \frac{\sigma_{b,\varepsilon}}{g_2^{-1}[g_2(\sigma_{b,\varepsilon}) + \varepsilon/(1-\varepsilon)]} - 1 \quad \text{and}$$

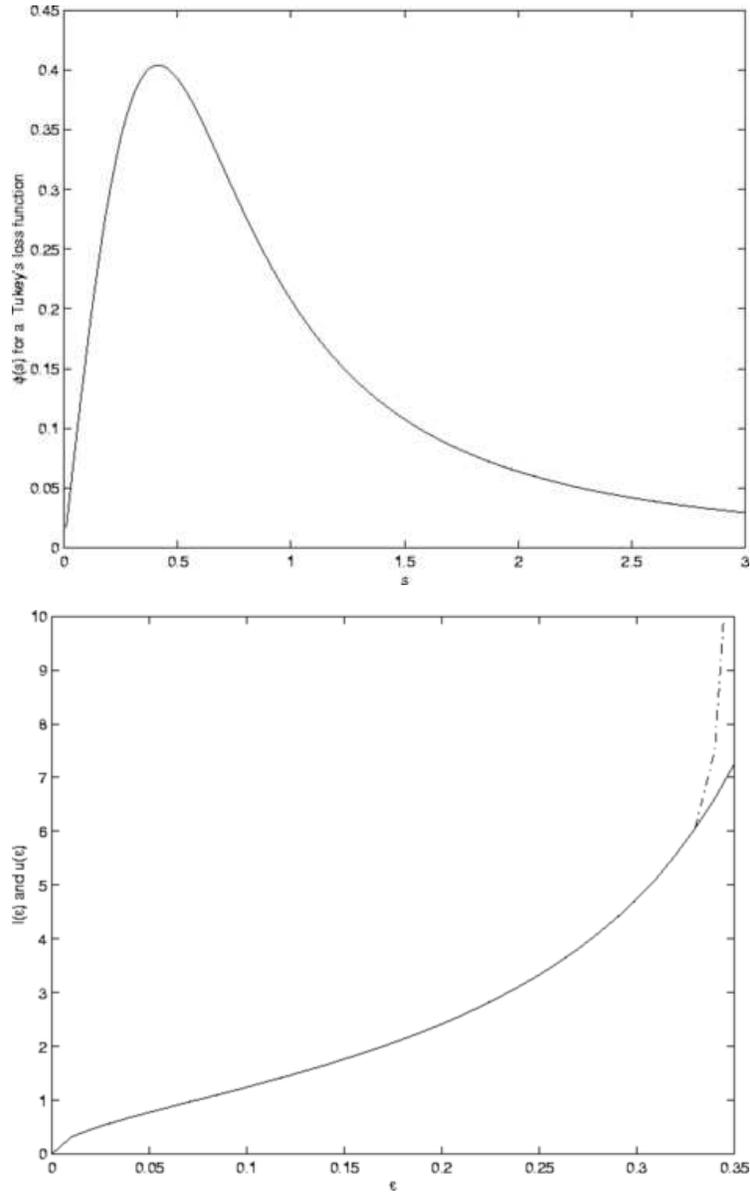

FIG. 1. *The graph on the top represents the function $\phi(s)$ for a biweight $\rho$-function. The graph on the bottom gives the maximum bias bounds [$\ell(\varepsilon)$, solid line; $u(\varepsilon)$, dotted-dashed line] for an MM-estimate based on biweight loss functions with 50% breakdown point and 95% efficiency under the Gaussian model.*

$$u(\varepsilon) = \frac{\gamma_{b,\varepsilon}}{g_2^{-1}[g_2(\gamma_{b,\varepsilon}) + \varepsilon/(1-\varepsilon)]} - 1.$$

The conditions given in (4.19), Theorem 4.1 and Theorem 4.2 for the Gaussian model are also being assumed here for (4.23), (4.24) and (4.25), respectively, for the Cauchy model. For an *MM*-estimate of regression, condition A4 can again be shown to hold when using a biweight loss function.

It is somewhat surprising that the expressions for $B_S(\varepsilon)$, $B_{CM}(\varepsilon)$ and $B_{MM}(\varepsilon)$ are of order $o(\varepsilon)$ as $\varepsilon \to 0$ under the Cauchy model, in contrast with the usual $\sqrt{\varepsilon}$ order. This is not, however, a contradiction of known results which establish general $\sqrt{\varepsilon}$ order for the maximum bias functions of regression estimates based on residuals since such results require either elliptical regressors, as in Yohai and Zamar [13], or the existence of second moments for the regressors, as in He [3] or Yohai and Zamar [14].

## 5. Maximum bias curve comparisons.

5.1. *The Gaussian model.* Most estimators need to be tuned so that they perform reasonably well under some important model, as well as being robust to deviations from the model. In practice, one often tunes an estimate so that it has good efficiency under the Gaussian model as well as a high breakdown point. For smooth $\rho$-functions, both the *MM*- and *CM*-estimates of regression can be tuned to have a 50% breakdown point and 95% asymptotic relative efficiency under the Gaussian model. This is also true for the class of $\tau$-estimates; see Yohai and Zamar [12] for details. Thus, these estimates cannot be ranked on the basis of their efficiency and breakdown point alone. Comparing their maximum bias behavior under the Gaussian model gives further insight into how these estimates are affected by deviations from the model.

Here, we again consider the estimates associated with the family of Tukey's biweight loss function (2.8). The 95% efficient biweight *MM*-estimate with a 50% breakdown point was discussed in the previous subsection. A 95% efficient biweight *CM*-estimate with a 50% breakdown point is obtained by choosing $\rho(u) \doteq \rho_T(u)$, $b = 0.5$ and the tuning constant $c = 4.835$; see [8] for details. In contrast, a 95% efficient biweight *S*-estimate of regression has a 12% breakdown point, whereas a biweight *S*-estimate with a 50% breakdown point is only 28.7% efficient under the Gaussian model.

Figure 2 represents the maximum bias functions under the Gaussian model of the *MM*-, *CM*- and $\tau$-estimates based on biweight functions and tuned so that they have 95% (asymptotic) efficiency under the Gaussian model and a 50% breakdown point, as well as the 95% efficient biweight *S*-estimate. We observe that up to $\varepsilon \approx 0.28$, the $\tau$-estimate has a larger bias than the *MM*-estimate and then a smaller bias afterward. The $\tau$-estimate, though, has a larger bias than the *CM*-estimate over essentially the entire range of $\varepsilon$. Up to $\varepsilon \approx 0.20$, *MM*- and *CM*-estimates are roughly equivalent, although for larger fractions of contamination, the *CM*-estimate is clearly better.

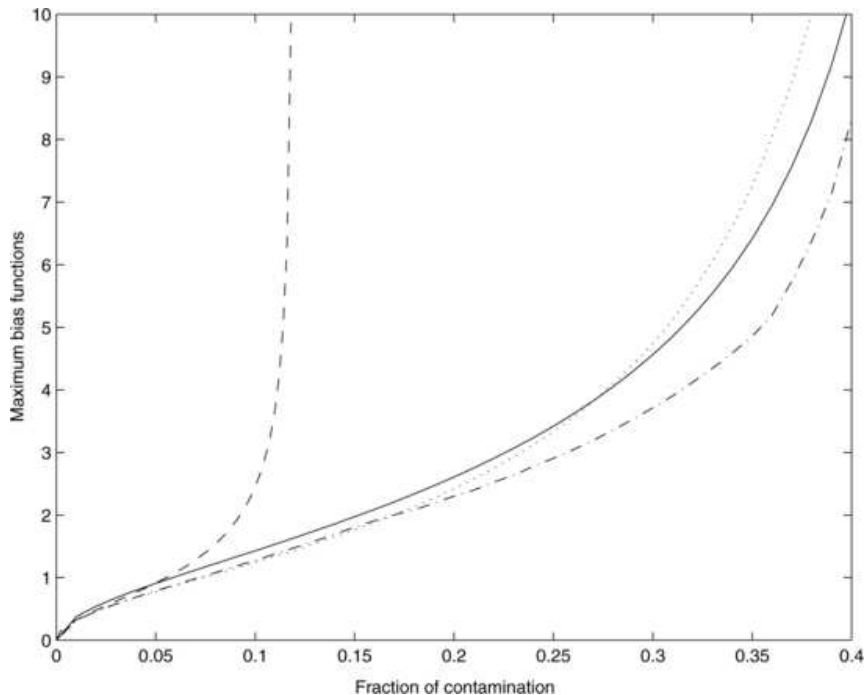

Fig. 2. *Maximum bias functions for a biweight S-estimate (dashed line), MM-estimate (dotted line, lower bound), $\tau$-estimate (solid line) and CM-estimate (dashed-dotted line). All of the estimates have 95% efficiency under the Gaussian model. The S-estimate has a breakdown point of* 12%, *whereas the others have a* 50% *breakdown point.*

As a further comparison, Figure 3 again shows the maximum bias function under the Gaussian model of the above 95% efficient biweight *MM*- and *CM*-estimates, as well as of the less efficient 50% breakdown point biweight *S*-estimate. Also included in Figure 3 is the biweight *CM*-estimate having a 50% breakdown point and an asymptotic relative efficiency of 61.1% under the Gaussian model, which corresponds to choosing the tuning constant $c = 2.568$. (The efficiency of the *CM*-estimate based on a biweight function with $b = 1/2$ and $c = 2.568$ under the Gaussian model is incorrectly reported as 28.7% rather than 61.1% in Table 1 of Mendes and Tyler [8]. The rest of Table 1 of [8] is correct.)

The maximum bias of the 95% efficient *MM*-estimate is uniformly larger than that of the corresponding *S*-estimate. This is consistent with the general result given in Theorem 3.1. The increase in bias for the *MM*-estimate is compensated by its increase in efficiency. A curious observation, though, is that for large fractions of contamination, the maximum bias of the 95% efficient $CM$-estimate is lower than that of the 28.7% efficient $S$-estimate. Furthermore, the maximum bias of the 61.1% efficient $CM$-estimate is almost identical to (and as shown theoretically in the next section, is never larger than) that of the 28.7% efficient $S$-estimate. That is, there is no trade-off between increased efficiency and maximum bias for this $CM$-estimate relative to the $S$-estimate. In practice, given that the maximum bias function of the 95% efficient $CM$-estimate does not greatly differ from that of the 61.1% estimate, the 95% efficient estimate would be preferable.

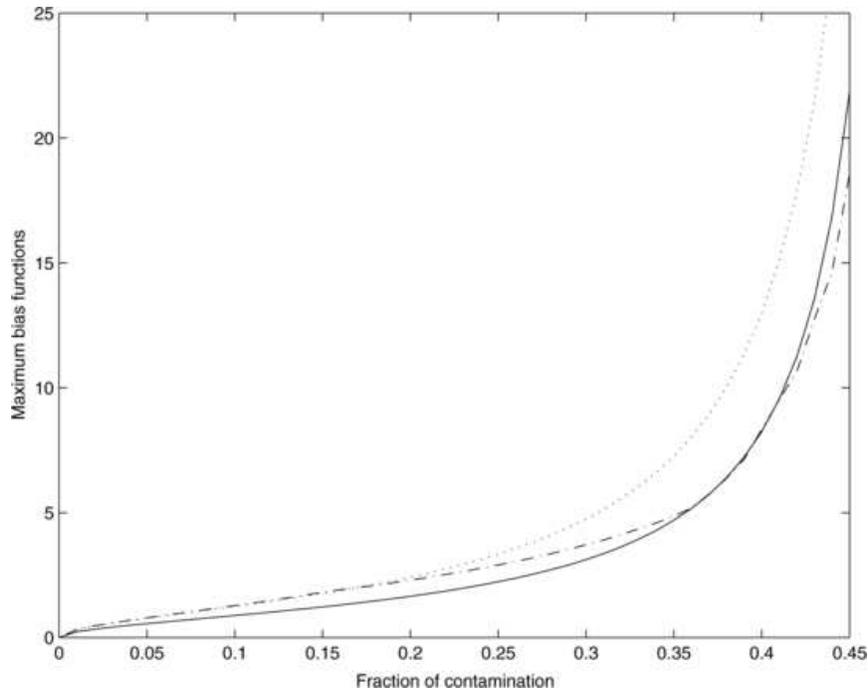

Fig. 3. *Maximum bias functions for a biweight S-estimate (solid line), MM-estimate (dotted line, lower bound), and two CM-estimates (dotted-dashed line and solid line). The plot for the S-estimate and the second CM-estimate are almost identical. All estimates have a* 50% *breakdown point. The MM-estimate and the first CM-estimates (dotted-dashed line) have* 95% *efficiency under the Gaussian model. The second CM-estimate (solid line) has an efficiency of* 61.1%, *whereas the efficiency of the S-estimate is* 28.7%.

5.2. *The Cauchy model.* We now consider the maximum bias behavior of $S$-, $MM$- and $CM$-estimates under the Cauchy model. Figure 4 shows the maximum bias functions under the Cauchy model for the $MM$- and $CM$-estimates which are 95% efficient under the Gaussian model, as well as for the 28.7% efficient biweight *S*-estimate and the 61.1% efficient *CM*-estimate discussed in Section 5.1. The breakdown point of each of these estimates remains 50% under the Cauchy model. The estimates, though, are not retuned here for the Cauchy model. Rather, our goal is to make further comparisons between the same estimates. In practice, given a specific estimate, one would wish to evaluate its robustness properties under various scenarios. From Figure 4, it can be noted that the 95% efficient *CM*-estimate tends to have the better maximum bias behavior under the Cauchy model, even better than that of the 61.1% efficient *CM*-estimate.

5.3. *Other considerations.* Aside from maximum bias functions, a classical way of evaluating the robustness of an estimate as it deviates from normality is to consider its efficiency under other distributions. The asymptotic efficiencies under the Gaussian model discussed in Section 5.1 depend on the distribution of the error term being normal. They do not however depend on the distribution of the carriers being normal. Rather, the only assumption needed for the carriers is that they possess second moments. This is also true for the asymptotic efficiencies at other symmetric error distributions; see, for example, Maronna, Bustos and Yohai [6]. In particular, the authors note that the asymptotic variance-covariance matrix of $\widehat{\boldsymbol{\theta}} = \mathbf{T}_n$ has the form $\sigma_u^2 \boldsymbol{\Sigma}_{\mathbf{x}}$, where $\boldsymbol{\Sigma}_{\mathbf{x}}$ is the variance-covariance matrix of the carriers $\mathbf{x}$ and $\sigma_u$ depends only on the distribution of the error term $u$.

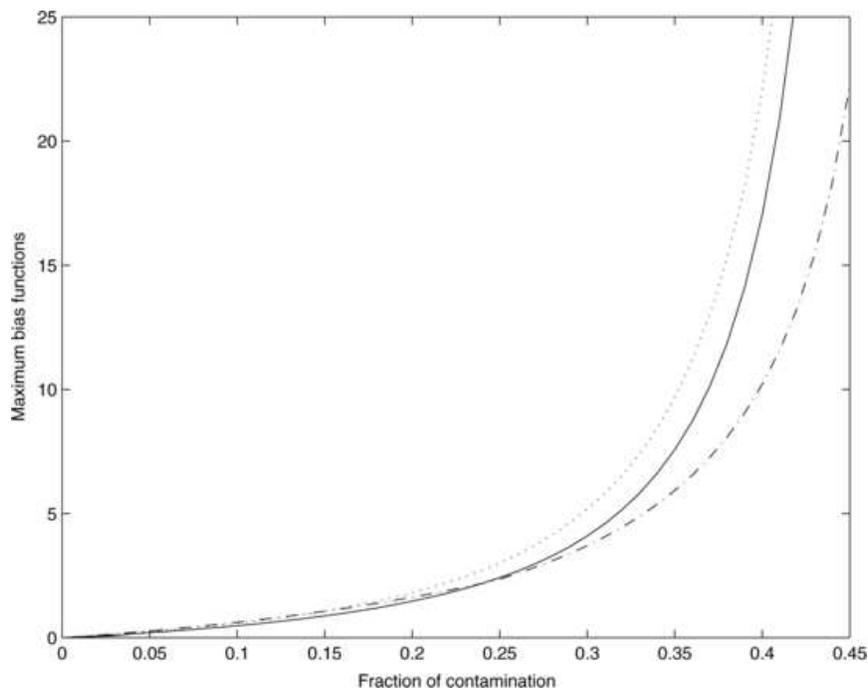

FIG. 4. *Maximum bias functions for a biweight S-estimate (solid line), an MM-estimate (dotted line, lower bound) and two CM-estimates (dotted-dashed line and solid line). The plot for the S-estimate and the second CM-estimate are almost identical. All estimates have a 50% breakdown point under the Cauchy model. The MM-estimate and the first CM-estimate have 95% efficiency under the Gaussian model, whereas the second CM-estimate and the S-estimate have efficiencies of 61.1% and 28.7%, respectively, under the Gaussian model.*

TABLE 1
*Asymptotic variances and residual gross error sensitivities of some S-, MM- and CM-estimates of regression under symmetric error distributions*

|      |      | NORM  | SL    | CAU   | T3    | DE    | CN    | UNIF    |
|------|------|-------|-------|-------|-------|-------|-------|---------|
| S95  | AVAR | 1.053 | 1.798 | 2.209 | 1.257 | 1.429 | 1.091 | 0.771   |
|      | RGES | 1.770 | 3.277 | 3.716 | 2.146 | 2.258 | 1.942 | 1.415   |
| MM95 | AVAR | 1.053 | 1.230 | 1.312 | 1.221 | 1.368 | 1.087 | 0.713   |
|      | RGES | 1.770 | 2.146 | 2.243 | 1.953 | 2.038 | 1.844 | 1.548   |
| CM95 | AVAR | 1.053 | 1.159 | 1.202 | 1.227 | 1.396 | 1.088 | 0.755   |
|      | RGES | 1.770 | 1.995 | 2.061 | 1.988 | 2.138 | 1.835 | 1.439   |
| CM61 | AVAR | 1.637 | 1.330 | 1.059 | 2.091 | 1.528 | 2.891 | 1.128   |
|      | RGES | 1.838 | 1.900 | 1.765 | 2.285 | 2.045 | 2.619 | 1.405   |
| S28  | AVAR | 3.484 | 1.330 | 1.059 | 2.091 | 1.528 | 2.891 | 120.336 |
|      | RGES | 2.850 | 1.900 | 1.765 | 2.285 | 2.045 | 2.619 | 15.621  |

In Table 1, we again consider the 95% efficient biweight *S*-, *MM*- and *CM*-estimates, the 28.7% efficient biweight *S*-estimate and the 61.1% efficient *CM*-estimate discussed in Section 5.1, where the efficiency is taken under a normal error model. These estimates are labeled *S95*, *MM95*, *CM95*, *S28* and *CM61*, respectively. For these estimates, we compute their asymptotic variances $\sigma_u^2$ (AVAR) under a variety of symmetric error models. Besides the standard normal (NORM), these models include the slash (SL), the Cauchy (CAU), the $t_3$-distribution (T3), the double-exponential (DE), a 90–10% mixture of a standard normal and a normal with mean zero and variance 9 (CN) and the uniform distribution on $(-1, 1)$ (UNIF). Each of these distributions is normalized so that its interquartile range is equal to that of the standard normal, namely 1.3490. This corresponds to multiplying the SL, CAU, T3, DE, CN or UNIF random variables by 0.4587, 0.6745, 0.8818, 0.9731, 0.9248 and 1.3490, respectively. Also included in Table 1 are the residual gross error sensitivities (RGES); see Hampel et al. [2]. Formulas for AVAR and RGES can be found in [8].

From Table 1, it can be noted that the estimates *MM95* and *CM95* behave similarly with respect to asymptotic variance and residual gross error sensitivity, with *CM95* being slightly better at the longer-tailed slash and Cauchy distributions and the *MM95* being slightly better at the more moderate $t_3$ and double-exponential distributions. Both *MM95* and *CM95* perform better than *S95* at longer-tailed distributions. The behavior of *S28* and *CM61* is the same, except at the normal and uniform distributions. At longer-tailed distributions, equality tends to hold for the constraint (2.9) on *CM61* and so as an estimate, it is asymptotically equivalent to *S28* under these distributions. Under the normal and the uniform distributions, there is a considerable difference in favor of *CM61*. Curiously, the behavior of *S28* and *CM61* under the Cauchy distribution is better than that of *MM95* and *CM95*. However, based on the overall behavior of the asymptotic variances and residual gross error sensitivities alone, either *MM95* and *CM95* would be preferable in practice.

**6. Bias-inadmissibility of *S*-estimates under the Gaussian model.** Throughout this section, the Gaussian model is presumed to hold, that is, it is presumed that both the response and the carriers are normal. In Section 5.1, it was noted that under the Gaussian model, the maximum bias function of the 61.1% efficient biweight *CM*-estimate is never smaller than that of the 27.78% efficient biweight *S*-estimate. In this section, we verify this result theoretically, rather than computationally. Moreover, we note that this result is not specific to the use of the biweight estimates. In general, we show that for a given *S*-estimate, it is usually possible to tune the corresponding *CM*-estimates (through the value of $c$) so that $B_{CM}(\varepsilon) \leq B_S(\varepsilon)$ for all $\varepsilon$ and with strict inequality for at least one value of $\varepsilon$. In such a case, we will say that with respect to the maximum bias criterion, the estimate $T_S$ is *inadmissible* under the Gaussian model since it can be *dominated* by $T_{CM}$.

To show this, we need to carefully compare the maximum bias functions of the *CM*-estimates and the *S*-estimates. An alternative representation for $B_{CM}(\varepsilon)$ in terms of $B_S(\varepsilon)$ under the normal model [see equations (4.21) and (4.19)] is given by

$$(6.26) \qquad \log[1 + B_{CM}^2(\varepsilon)] = \log[1 + B_S^2(\varepsilon)] + 2d_c(\varepsilon),$$

where $d_c(\varepsilon) = h_c(\varepsilon, \gamma_{b,\varepsilon}) - h_c(\varepsilon, \sigma_{b,\varepsilon})$ and

$$(6.27) \qquad h_c(\varepsilon, \sigma) = A_{c,\varepsilon}(\sigma) - \inf_{s \geq \sigma} A_{c,\varepsilon}(s).$$

The functionals $T_{CM}$ and $T_S$ in (6.26) are understood to be defined by using the same $\rho$ and the same value of $b$. From representation (6.26), it is apparent that what we need to consider is the sign of $d_c(\varepsilon)$ in terms of $c$ and $\varepsilon$. The following result represents a first step in determining appropriate values of the tuning constant $c$ necessary for showing the bias inadmissibility of an *S*-estimate. The value of $K$ below is defined within condition A4.

THEOREM 6.1. *Suppose that $\rho$ is such that conditions* A3 *and* A4 *hold. Then:*

(i) *if $c \leq 1/K$, then $B_{CM}(\varepsilon) = B_S(\varepsilon)$ for all $\varepsilon$;*

(ii) *for any $\varepsilon$ such that $c > c(\varepsilon) \doteq \varepsilon^{-1} \log(\sigma_{b,\varepsilon}/\gamma_{b,\varepsilon})$, we have $B_{CM}(\varepsilon) > B_S(\varepsilon)$.*

As a consequence, for the $CM$-estimate to improve upon the maximum bias function of the $S$-estimate, one needs to choose $c > 1/K$. On the other hand, if $c_o \doteq \inf\{c(\varepsilon): 0 < \varepsilon < b\}$, then we also need to choose $c \leq c_o$. This range is not empty since, as shown in the Appendix,

$$(6.28) \qquad 1/K < (1-b)/K + b/\phi(\gamma_{b,0}) \leq c_o,$$

where $\phi(s)$ is defined within condition A4.

For $c \leq 1/K$, the $CM$-functional is the same as the $S$-functional at $H_o$, as well as at any $H$ in an $\varepsilon$-contaminated neighborhood of $H_o$. This is because equality is obtained in the constraint (2.9) for the $CM$-estimate and when equality is obtained, the $CM$-estimate gives the same solution as the corresponding $S$-estimate. Thus, for $c \leq 1/K$, the $CM$-estimate has the same maximum bias function as the corresponding $S$-estimate. On the other hand, for large values of $c$, the $CM$-estimate tends to give a solution similar to the least squares solution, so one expects the maximum bias function to be unacceptably large, even though the breakdown point may be close to $1/2$. In fact, one can note from (4.21) that for any $\varepsilon$, $B_{CM}(\varepsilon) \to \infty$ as $c \to \infty$.

Varying the tuning constant $c$ may decrease the maximum bias for some values of $\varepsilon$, while increasing the maximum bias for other values of $\varepsilon$. The question we address now is whether it is possible to find a moderate value of $c$ (necessarily between $1/K$ and $c_o$) such that the maximum bias function of the $CM$-estimate improves upon the maximum bias function of the $S$-estimate.

The following result shows that in most cases of interest, the condition $c \leq c_o$ is not only necessary, but also sufficient, to obtain $B_{CM}(\varepsilon) \leq B_S(\varepsilon)$ for all $\varepsilon$. The value of $\sigma_M$ below is also defined within condition A4.

THEOREM 6.2. *Suppose that the assumptions of Theorem* 6.1 *hold. If $c \leq c_o$ and $g(\sigma_M) \leq b$, then $B_{CM}(\varepsilon) \leq B_S(\varepsilon)$ for all $\varepsilon > 0$.*

REMARK 6.1. This result cannot be improved upon. That is, if $c > c(\varepsilon)$, then $B_S(\varepsilon) < B_{CM}(\varepsilon)$ by Theorem 6.1. Also, if $c \leq c_o$ and $g(\sigma_M) > b$, then either $B_S(\varepsilon) < B_{CM}(\varepsilon)$ for some $\varepsilon$ or $B_S(\varepsilon) = B_{CM}(\varepsilon)$ for all $\varepsilon$. This remark is verified in the Appendix.

In order to show that an *S*-estimate can be *dominated* by a *CM*-estimate with $c$ chosen so that $1/K < c \leq c_o$, it remains to be shown that for some $0 < \varepsilon < b$, $B_{CM}(\varepsilon) < B_S(\varepsilon)$. For specific examples, this can be checked numerically; under additional assumptions, though, it can be shown analytically.

THEOREM 6.3. *Suppose that the assumptions of Theorem 6.2 hold. Furthermore, suppose that $g(s)$ is convex and*

$$\phi(\sigma_{b,0}) \geq \frac{[1 - g(\sigma_M)]^2 (1 - b)}{2 - [b + g(\sigma_M)]}. \tag{6.29}$$

*Then for any value $c$ such that*

$$c_1 \doteq \frac{\log(\sigma_M/\sigma_{b,0})}{b - g(\sigma_M)} < c \leq \frac{1}{\phi(\sigma_{b,0})} = c_o,$$

*the CM-estimate of regression* **dominates** *the S-estimate of regression with respect to the maximum bias function. Furthermore, this range of values for $c$ is not empty.*

REMARK 6.2. From the proof of Theorem 6.3, it follows that a condition more general than (6.29) under which the conclusions also hold is $c_o = \lim_{\varepsilon \to 0^+} c(\varepsilon)$. However, (6.29) is easier to check and holds in most cases of interest.

Consider the biweight *S*-estimate with breakdown point $b \leq 1/2$. It can be verified that the conditions of Theorem 6.3 hold whenever $b > 0.410$, so any such biweight *S*-estimate is inadmissible with respect to maximum bias under the Gaussian model. For $b = 1/2$, that is, for the 27.78% efficient biweight *S*-estimate, the value of $c = 2.568$ falls within the interval given in Theorem 6.3. Hence, the 61.1% efficient biweight *CM*-estimate *dominates* the 27.78% efficient biweight *S*-estimate with respect to maximum bias under the Gaussian model. As noted in Section 5.1, although the decrease in maximum bias is negligible, the increase in efficiency is not.

As another example, consider the $\alpha$-quantile regression estimates. These correspond to *S*-estimates with $\rho(u) = I\{|u| \geq 1\}$ and $b = 1 - \alpha$. It is straightforward to verify that the conditions of Theorem 6.3 hold in this case whenever $b > 0.3173$, so the $\alpha$-quantile regression estimates with $\alpha < 0.6837$ are inadmissible under the Gaussian model with respect to maximum bias. Again, the decrease in maxbias is not large. For example, in the special case $\alpha = b = 0.5$, for which the resulting $\alpha$-quantile estimate corresponds to Rousseeuw's [9] least median of squares estimate (*LMS*), the best improvement is only 95.7% of the *LMS* bias.

The $\alpha$-quantile estimates are often referred to as minimax bias regression estimates. Martin, Yohai and Zamar [7] show that within the class of $M$-estimates of regression with general scale, an $\alpha$-quantile estimate minimizes the maximum bias at $\varepsilon$, with the value of $\alpha$ depending on $\varepsilon$. Yohai and Zamar [13] generalize this minimax result to the class of all residual admissible estimates of regression. Under the Gaussian model, an $\alpha$-quantile estimate can be shown to have minimax bias for some $\varepsilon$ whenever $0.500 < \alpha < 0.6837$ or, equivalently, whenever $0.3173 < b < 0.500$. Despite having minimax bias under the Gaussian model for a given $\varepsilon$, these $\alpha$-quantile regression estimates are still inadmissible under the Gaussian model with respect to maximum bias. In particular, as shown in Theorem 6.1, for a given $\alpha$-quantile estimate having minimax bias at $\varepsilon = \varepsilon_\alpha$, it is possible to construct a $CM$-estimate which also has the same maximum bias at $\varepsilon = \varepsilon_\alpha$. Moreover, the maximum bias of this $CM$-estimate is never larger than the maximum bias of the given $\alpha$-quantile estimate at any other value of $\varepsilon \neq \varepsilon_\alpha$, and furthermore is smaller than the maximum bias of the given $\alpha$-quantile estimate at some values of $\varepsilon \neq \varepsilon_\alpha$. Although the decrease in the maximum bias may not be of practical importance, these observations expose some limitations of the notion of minimax bias.

The minimax bias results given in [13] for the $\alpha$-quantile regression estimates apply more generally than to just the Gaussian model. They also apply to models having a symmetric unimodal error term along with elliptically distributed carriers. Under such models, though, we conjecture that the $\alpha$-quantile regression estimates may again be inadmissible with respect to maximum bias, but we do not pursue this topic further here. The value of $\alpha$ which attains the minimum maxbias at a particular $\varepsilon$ is not only dependent on the value of $\varepsilon$, but also dependent on the particular model. That is, a particular $\alpha$-quantile estimate is not necessarily minimax at $\varepsilon$ over a range of models, but is only known to be minimax at $\varepsilon$ under a specific model. Any estimate which can be shown to *dominate* an $\alpha$-quantile estimate would most likely need to be model-specific.

## APPENDIX

In this section, we include proofs of the results and other technical questions.

PROOF OF THEOREM 3.1. It can be shown, following the proofs of Lemmas 4, 5 and 6 in [1], that for all $s > 0$ and $t \in \mathbb{R}$, there exist $\alpha_t \in \mathbb{R}$ and $\boldsymbol{\theta}_t \in \mathbb{R}^p$ such that

$$m(t,s) = \mathrm{E}_{H_o}\rho_2\left(\frac{y - \alpha_t - \mathbf{x}'\boldsymbol{\theta}_t}{s}\right) - \mathrm{E}_{H_o}\rho_2\left(\frac{y}{s}\right).$$

Also, we can show that $m(t, s)$ is a strictly increasing function of $t$ for all $s > 0$. It follows that $h_1(t)$ is also strictly increasing.

We first show that $B_{MM}(\varepsilon) \leq t_2$, where $t_2$ is such that $h_2(t_2) = \varepsilon/(1-\varepsilon)$. Let $\tilde{\boldsymbol{\theta}} \in \mathbb{R}^p$ be such that $\tilde{t} = \|\tilde{\boldsymbol{\theta}}\| > t_2$. We shall prove that

$$\text{E}_H \rho_2 \left( \frac{y - \alpha - \mathbf{x}'\tilde{\boldsymbol{\theta}}}{s(H)} \right) > \text{E}_H \rho_2 \left( \frac{y}{s(H)} \right) \quad \text{for each } \alpha \in \mathbb{R} \text{ and } H \in V\varepsilon. \tag{A.30}$$

Let $H = (1-\varepsilon)H_o + \varepsilon \tilde{H}$. We have that

$$m[\tilde{t}, s(H)] > m[t_2, s(H)] \geq \inf_{\underline{s} \leq s \leq \overline{s}} m(t_2, s) = h_2(t_2) = \frac{\varepsilon}{1-\varepsilon}.$$

Therefore, for each $\alpha \in \mathbb{R}$ and $H \in V\varepsilon$,

$$\text{E}_{H_o} \rho_2 \left( \frac{y - \alpha - \mathbf{x}'\tilde{\boldsymbol{\theta}}}{s(H)} \right) - \text{E}_{H_o} \rho_2 \left( \frac{y}{s(H)} \right) > \frac{\varepsilon}{1-\varepsilon},$$

that is,

$$(1-\varepsilon) \text{E}_{H_o} \rho_2 \left( \frac{y - \alpha - \mathbf{x}'\tilde{\boldsymbol{\theta}}}{s(H)} \right) > (1-\varepsilon) \text{E}_{H_o} \rho_2 \left( \frac{y}{s(H)} \right) + \varepsilon.$$

It follows that for every $\alpha \in \mathbb{R}$ and $H \in V_\varepsilon$,

$$\text{E}_H \rho_2 \left( \frac{y - \alpha - \mathbf{x}'\tilde{\boldsymbol{\theta}}}{s(H)} \right) \geq (1-\varepsilon) \text{E}_{H_o} \rho_2 \left( \frac{y - \alpha - \mathbf{x}'\tilde{\boldsymbol{\theta}}}{s(H)} \right)$$

$$> (1-\varepsilon) \text{E}_{H_o} \rho_2 \left( \frac{y}{s(H)} \right) + \varepsilon \geq \text{E}_H \rho_2 \left( \frac{y}{s(H)} \right),$$

that is, inequality (A.30) holds. The last inequality above follows from A3(ii).

Next, we show that $B_{MM}(\varepsilon) \geq t_1$, where $t_1$ is such that $h_1(t_1) = \varepsilon/(1-\varepsilon)$. Since $B_S(\varepsilon) < t_1$, we can select an arbitrary $t > 0$ such that $B_S(\varepsilon) < t < t_1$. It is enough to show that $B_{MM}(\varepsilon) \geq t$. We know that there exist $\alpha_t \in \mathbb{R}$ and $\boldsymbol{\theta}_t \in \mathbb{R}^p$ such that

$$h_1(t) = m(t, \overline{s}) = \text{E}_{H_o} \rho_2 \left( \frac{y - \alpha_t - \mathbf{x}'\boldsymbol{\theta}_t}{\overline{s}} \right) - \text{E}_{H_o} \rho_2 \left( \frac{y}{\overline{s}} \right).$$

Since $h_1$ is strictly increasing, $h_1(t) < h_1(t_1) = \varepsilon/(1-\varepsilon)$. It follows that

$$(1-\varepsilon) \text{E}_{H_o} \rho_2 \left( \frac{y - \alpha_t - \mathbf{x}'\boldsymbol{\theta}_t}{\overline{s}} \right) < (1-\varepsilon) \text{E}_{H_o} \rho_2 \left( \frac{y}{\overline{s}} \right) + \varepsilon. \tag{A.31}$$

Define the following sequence of contaminating distributions: $\tilde{H}_n = \delta_{(y_n, \mathbf{x}_n)}$, where $\mathbf{x}_n = n\boldsymbol{\theta}_t$ and $y_n = \alpha_t + \mathbf{x}'_n \boldsymbol{\theta}_t = \alpha_t + nt^2$. Let $H_n = (1-\varepsilon)H_o + \varepsilon \tilde{H}_n$ and $\boldsymbol{\theta}_n = \mathbf{T}(H_n)$. Suppose that $\sup_n \|\boldsymbol{\theta}_n\| < t$, in order to produce a contradiction. Under this assumption, there exists a convergent subsequence, denoted also by $\{\boldsymbol{\theta}_n\}$, such that $\lim_{n\to\infty} \boldsymbol{\theta}_n = \tilde{\boldsymbol{\theta}}$, where $\|\tilde{\boldsymbol{\theta}}\| = \tilde{t} < t$. Assume for a moment that the sequence of intercept functionals evaluated at $H_n$, $\alpha_n = t(H_n)$, satisfies $\lim_{n\to\infty} |\alpha_n| = \infty$. Then

$$\lim_{n\to\infty} \mathrm{E}_{H_n} \rho_2\left(\frac{y - \alpha_n - \mathbf{x}'\boldsymbol{\theta}_n}{s(H_n)}\right) = (1-\varepsilon) + \varepsilon \lim_{n\to\infty} \rho_2\left(\frac{y_n - \alpha_n - \mathbf{x}'_n\boldsymbol{\theta}_n}{s(H_n)}\right)$$

$$> (1-\varepsilon) \lim_{n\to\infty} \mathrm{E}_{H_o} \rho_2\left(\frac{y - \alpha_t - \mathbf{x}'\boldsymbol{\theta}_t}{s(H_n)}\right)$$

$$= \lim_{n\to\infty} \mathrm{E}_{H_n} \rho_2\left(\frac{y - \alpha_t - \mathbf{x}'\boldsymbol{\theta}_t}{s(H_n)}\right),$$

but this fact contradicts the definition of $(\alpha_n, \boldsymbol{\theta}_n)$. Note that $0 < \underline{s} < s(H_n) < \overline{s} < \infty$ implies $\lim_{n\to\infty} \mathrm{E}_{H_o} \rho_2[(y - \alpha_t - x'\boldsymbol{\theta}_t)/s(H_n)] < 1$ which, in turn, implies the strict inequality above. Therefore, we can assume, without loss of generality, that $\lim_{n\to\infty} \alpha_n = \tilde{\alpha}$ for some finite $\tilde{\alpha} \in \mathbb{R}$. As a consequence, we have

(A.32)
$$\lim_{n\to\infty} \left|\frac{y_n - \alpha_n - \mathbf{x}'_n\boldsymbol{\theta}_n}{s(H_n)}\right| = \infty \quad \text{and}$$

$$\left|\frac{y_n - \alpha_t - \mathbf{x}'_n\boldsymbol{\theta}_t}{s(H_n)}\right| = 0 \quad \text{for each } n.$$

We now prove that $\lim_{n\to\infty} s(H_n) = \overline{s}$ for any convergent subsequence $s(H_n)$. Let $s_\infty = \lim_{n\to\infty} s(H_n)$. Note that $\overline{s}$ satisfies the equation

(A.33) $$(1-\varepsilon)\mathrm{E}_{H_o}\rho_1(y/\overline{s}) + \varepsilon = b.$$

Let $(\gamma_n, \boldsymbol{\beta}_n) = (t_1(H_n), \mathbf{T}_1(H_n))$ be the regression $S$-estimate based on $\rho_1$. We know that $\|\boldsymbol{\beta}_n\| \leq B_S(\varepsilon) < t$ for all $n$, so that, without loss of generality, $\lim_{n\to\infty} \boldsymbol{\beta}_n = \tilde{\boldsymbol{\beta}}$, where $\|\tilde{\boldsymbol{\beta}}\| < t$. Assume that $\lim_{n\to\infty} |\gamma_n| = \infty$. Since

(A.34) $$\mathrm{E}_{H_n}\rho_1\left(\frac{y - \gamma_n - \mathbf{x}'\boldsymbol{\beta}_n}{s(H_n)}\right) = b,$$

letting $n \to \infty$, it follows that

$$b = \lim_{n\to\infty} \mathrm{E}_{H_n}\rho_1\left(\frac{y - \gamma_n - \mathbf{x}'\boldsymbol{\beta}_n}{s(H_n)}\right) = (1-\varepsilon) + \varepsilon \lim_{n\to\infty} \rho_1\left(\frac{y_n - \gamma_n - \mathbf{x}'_n\boldsymbol{\beta}_n}{s(H_n)}\right)$$

$$> (1-\varepsilon) \lim_{n\to\infty} \mathrm{E}_{H_o}\rho_1\left(\frac{y - \alpha_t - \mathbf{x}'\boldsymbol{\theta}_t}{s(H_n)}\right) = \mathrm{E}_{H_n}\rho_1\left(\frac{y - \alpha_t - \mathbf{x}'\boldsymbol{\theta}_t}{s(H_n)}\right).$$

Then there exists $s_n < s(H_n)$ such that

$$\mathrm{E}_{H_n}\rho_1\left(\frac{y - \alpha_t - \mathbf{x}'\boldsymbol{\theta}_t}{s_n}\right) = b,$$

but this fact contradicts the definition of $(\gamma_n, \boldsymbol{\beta}_n)$. Therefore, we can also assume, without loss of generality, that $\lim_{n\to\infty} \gamma_n = \tilde{\gamma}$ for some finite $\tilde{\gamma} \in \mathbb{R}$. As a consequence, letting $n \to \infty$ in (A.34), we obtain

$$b = (1-\varepsilon)\mathrm{E}_{H_o}\rho_1\left(\frac{y - \tilde{\gamma} - \mathbf{x}'\tilde{\boldsymbol{\beta}}}{s_\infty}\right) + \varepsilon \geq (1-\varepsilon)\mathrm{E}_{H_o}\rho_1(y/s_\infty) + \varepsilon.$$

Comparing the last equation with (A.33), we deduce that $s_\infty \geq \overline{s}$. Since $\overline{s} = \sup_{H \in V_\varepsilon} s(H)$, we have $s_\infty = \overline{s}$. We use this fact to obtain Equations (A.35) and (A.36) below.

Equations (A.31) and (A.32) imply that

$$\lim_{n\to\infty} E_{H_n}\rho_2\left(\frac{y - \alpha_n - \mathbf{x}'\boldsymbol{\theta}_n}{s(H_n)}\right) = (1-\varepsilon)E_{H_o}\rho_2\left(\frac{y - \tilde{\alpha} - \mathbf{x}'\tilde{\boldsymbol{\theta}}}{\overline{s}}\right) + \varepsilon$$

(A.35)
$$\geq (1-\varepsilon)E_{H_o}\rho_2\left(\frac{y}{\overline{s}}\right) + \varepsilon$$

$$> (1-\varepsilon)E_{H_o}\rho_2\left(\frac{y - \alpha_t - \mathbf{x}'\boldsymbol{\theta}_t}{\overline{s}}\right).$$

On the other hand, applying (A.32), we have

(A.36) $$\lim_{n\to\infty} E_{H_n}\rho_2\left(\frac{y - \alpha_t - \mathbf{x}'\boldsymbol{\theta}_t}{s(H_n)}\right) = (1-\varepsilon)E_{H_o}\rho_2\left(\frac{y - \alpha_t - \mathbf{x}'\boldsymbol{\theta}_t}{\overline{s}}\right).$$

Therefore, for sufficiency large $n$,

$$E_{H_n}\rho_2\left(\frac{y - \alpha_n - \mathbf{x}'\boldsymbol{\theta}_n}{s(H_n)}\right) > E_{H_n}\rho_2\left(\frac{y - \alpha_t - \mathbf{x}'\boldsymbol{\theta}_t}{s(H_n)}\right).$$

This last inequality contradicts the definition of $(\alpha_n, \boldsymbol{\theta}_n)$. For every $t > 0$ such that $B_S(\varepsilon) < t < t_1$, we have found a sequence of distributions $\{H_n\}$ in the neighborhood $V_\varepsilon$ such that $\sup_n \|\mathbf{T}(H_n)\| \geq t$. Therefore, $B_{MM}(\varepsilon) \geq t_1$. □

PROOF OF THEOREM 3.2. It is enough to check that the functional $J(F)$ defined in (3.15) satisfies condition A1 in [1]. For instance, the monotonicity condition A1(a) follows immediately from the monotonicity of the $M$-scale $\sigma(F)$. □

PROOF OF THEOREM 4.1. We will apply Theorem 3.2. Let $F_o = (1-\varepsilon)F_{H_o,0,\mathbf{0}} + \varepsilon\delta_\infty$. Then $\sigma(F_o) = \sigma_{b,\varepsilon}$ and

(A.37) $$r_{CM}(\varepsilon) = J_{CM}(F_o) = \inf_{s \geq \sigma_{b,\varepsilon}} A_{c,\varepsilon}(s) + c\varepsilon.$$

On the other hand, if $\|\boldsymbol{\theta}\| = t$ and $F_t = (1-\varepsilon)F_{H_o,0,\boldsymbol{\theta}} + \varepsilon\delta_o$, then, when $H_o$ is multivariate normal, we have $\sigma(F_t) = (1+t^2)^{1/2}\gamma_{b,\varepsilon}$ and

$$(A.38) \quad m_{CM}(t) = J_{CM}(F_t) = \frac{1}{2}\log(1+t^2) + \inf_{s \geq \gamma_{b,\varepsilon}} A_{c,\varepsilon}(s).$$

From (3.18), we know that $B_{CM}(\varepsilon) = t_\varepsilon$, where $m_{CM}(t_\varepsilon) = r_{CM}(\varepsilon)$. Matching the expressions in Equations (A.37) and (A.38) and solving for $t$ yields the result. □

PROOF OF THEOREM 4.2. Let $t \in \mathbb{R}$ be arbitrary. Under the assumptions, the function $m(t,s)$ is continuously differentiable with respect to $s$, with derivative given by

$$\frac{\partial m(t,s)}{\partial s} = \frac{1}{s}\left[\phi_2(s) - \phi_2\left(\frac{s}{(1+t^2)^{1/2}}\right)\right].$$

Since $\phi_2(s)$ is unimodal, for each $\varepsilon$ we have that $m(t,s)$ is (a) strictly increasing for $s \in [\underline{s}, \overline{s}]$, (b) strictly decreasing for $s \in [\underline{s}, \overline{s}]$ or (c) it has a unique critical point $\tilde{s} \in (\underline{s}, \overline{s})$ which is a local maximum. In any of the three cases, the global minimum of $m(t,s)$ for $s \in [\underline{s}, \overline{s}]$ is attained at one of the two extremes of the interval. That is,

$$h_2(t) = \inf_{\underline{s} \leq s \leq \overline{s}} m(t,s) = \min\{m(t,\underline{s}), m(t,\overline{s})\}.$$

From Theorem 3.1, an upper bound for the maximum bias is given by the value of $t_\varepsilon$ such that $h_2(t_\varepsilon) = \varepsilon/(1-\varepsilon)$. If $h_2(t_\varepsilon) = m(t_\varepsilon, \overline{s})$, then $h_1(t_\varepsilon) = h_2(t_\varepsilon)$ and therefore $t_\varepsilon = \ell(\varepsilon)$. On the other hand, if $h_2(t_\varepsilon) = m(t_\varepsilon, \underline{s})$, then we have that $t_\varepsilon = u(\varepsilon)$. Hence, the result follows. □

PROOF OF (4.23). We apply Theorem 1 from [1]. Following the notation in that paper, we have that $c = \sigma_{b,\varepsilon}$. On the other hand,

$$m(t) \doteq \inf_{\|\boldsymbol{\theta}\|=t} \inf_{\alpha \in \mathbb{R}} J_S[(1-\varepsilon)F_{H_o,\alpha,\boldsymbol{\theta}} + \varepsilon\delta_0]$$
$$= \inf_{\|\boldsymbol{\theta}\|=t} J_S[(1-\varepsilon)F_{H_o,0,\boldsymbol{\theta}} + \varepsilon\delta_0] = \inf_{\|\boldsymbol{\theta}\|=t} S(\boldsymbol{\theta}),$$

where $S(\boldsymbol{\theta})$ is such that

$$(A.39) \quad (1-\varepsilon)\mathrm{E}_{H_0}\rho\left(\frac{y - \mathbf{x}'\boldsymbol{\theta}}{S(\boldsymbol{\theta})}\right) = b.$$

Since $y - \mathbf{x}'\boldsymbol{\theta}$ is distributed as $(1 + \sum_i |\theta_i|)Z$, where $Z$ is standard Cauchy, we have that (A.39) amounts to

$$(1-\varepsilon)g\left(\frac{S(\boldsymbol{\theta})}{1 + \sum_i |\theta_i|}\right) = b.$$

Therefore,

$$S(\boldsymbol{\theta}) = \left[1 + \sum_i |\theta_i|\right]\gamma_{b,\varepsilon}$$

and

(A.40) $$m(t) = \inf_{\|\boldsymbol{\theta}\|=t}\left[1 + \sum_i |\theta_i|\right]\gamma_{b,\varepsilon} = (1+t)\gamma_{b,\varepsilon}.$$

Finally, since $B_s(\varepsilon) = t$, where $m(t) = \sigma_{b,\varepsilon}$, the result follows from (A.40). □

PROOF OF (4.24). Clearly, under the Cauchy model, the expression for $r_{CM}(\varepsilon)$ is formally the same as that corresponding to the Gaussian model. We just have to compute $g(s)$ with respect to the Cauchy distribution instead of the normal. On the other hand, under the Cauchy model it is not difficult to check that

$$m_{CM}(t) = \log(1+t) + \inf_{s \geq \gamma_{b,\varepsilon}} A_{c,\varepsilon}(s),$$

where $A_{c,\varepsilon}(s)$ is defined by (4.20). Since the bias satisfies $m[B_{CM}(\varepsilon)] = r_{CM}(\varepsilon)$, the result follows. □

PROOF OF (4.25). The same arguments as above yield the following expression for the function $m(t,s)$ under the Cauchy model:

$$m(t,s) = g_2\left(\frac{s}{1+t}\right) - g_2(s).$$

From this expression, the computation of $\ell(\varepsilon)$ and $u(\varepsilon)$ under the Cauchy model is straightforward, as follows:

$$\ell(\varepsilon) = h_1^{-1}\left(\frac{\varepsilon}{1-\varepsilon}\right) = \frac{\sigma_{b,\varepsilon}}{g_2^{-1}[g_2(\sigma_{b,\varepsilon}) + \varepsilon/(1-\varepsilon)]} - 1$$

and

$$u(\varepsilon) = \frac{\gamma_{b,\varepsilon}}{g_2^{-1}[g_2(\gamma_{b,\varepsilon}) + \varepsilon/(1-\varepsilon)]} - 1.$$

Since we are assuming that $\phi(s)$ is unimodal, the same proof as in the case of the Gaussian model yields (4.25). □

PROOF OF THEOREM 6.1. (i) Computing the derivative of $A_{c,\varepsilon}(s)$ with respect to $s$, we see that $A_{c,\varepsilon}(s)$ is nondecreasing when $c < [(1-\varepsilon)\phi(s)]^{-1}$. Since $K^{-1} < [(1-\varepsilon)\phi(s)]^{-1}$ for all $\varepsilon$ and $s > 0$, the condition $c \leq K^{-1}$ implies that $A_{c,\varepsilon}(s)$ is nondecreasing for all $\varepsilon$ and $s > 0$. As a consequence, $h_c(\varepsilon, \sigma) = $ 0 for all $\varepsilon$ and $\sigma > 0$. Then $d_c(\varepsilon) = 0$ for all $\varepsilon$, which implies that $B_{CM}(\varepsilon) = B_S(\varepsilon)$ for all $\varepsilon$.

(ii) Since $g(\sigma_{b,\varepsilon}) = (b - \varepsilon)/(1 - \varepsilon)$ and $g(\gamma_{b,\varepsilon}) = b/(1 - \varepsilon)$, it follows that

(A.41) $$A_{c,\varepsilon}(\sigma_{b,\varepsilon}) < A_{c,\varepsilon}(\gamma_{b,\varepsilon}) \iff c > c(\varepsilon).$$

However, if $A_{c,\varepsilon}(\sigma_{b,\varepsilon}) < A_{c,\varepsilon}(\gamma_{b,\varepsilon})$, then $d_c(\varepsilon) > 0$ and hence $B_S(\varepsilon) < B_{CM}(\varepsilon)$. □

PROOF OF (6.28). By using implicit differentiation, one obtains

$$\frac{\partial \sigma_{b,\varepsilon}}{\partial \varepsilon} = \frac{1}{(1-\varepsilon)^2} \frac{(1-b)\sigma_{b,\varepsilon}}{\phi(\sigma_{b,\varepsilon})} \quad \text{and} \quad \frac{\partial \gamma_{b,\varepsilon}}{\partial \varepsilon} = \frac{1}{(1-\varepsilon)^2} \frac{-b\gamma_{b,\varepsilon}}{\phi(\gamma_{b,\varepsilon})}.$$

This then gives

(A.42) $$\frac{\partial \varepsilon c(\varepsilon)}{\partial \varepsilon} = \frac{1}{(1-\varepsilon)^2} \left( \frac{1-b}{\phi(\sigma_{b,\varepsilon})} + \frac{b}{\phi(\gamma_{b,\varepsilon})} \right) \geq \frac{1-b}{K} + \frac{b}{\phi(\sigma_{b,0})}.$$

The last inequality follows since, as noted previously, $\gamma_{b,\varepsilon} < \sigma_{b,0} < \sigma_M$. This then implies (6.28). □

PROOF OF THEOREM 6.2. For $c \leq 1/K$, it has already been noted that the maximum bias functions are the same, so we need only consider $1/K < c \leq c_o$. In general, for $c > 1/K$ and under Assumption A4, the function $A_{c,\varepsilon}(s)$ has the following properties:

(i) $A_{c,\varepsilon}(0) = -\infty$ and $A_{c,\varepsilon}(\infty) = \infty$;
(ii) $A_{c,\varepsilon}(s)$ has two critical points, say $\sigma_L(c, \varepsilon) \leq \sigma_U(c, \varepsilon)$, with
   $A_{c,\varepsilon}(s)\Uparrow$ over 0 to $\sigma_L(c, \varepsilon)$,
   $A_{c,\varepsilon}(s)\Downarrow$ over $\sigma_L(c, \varepsilon)$ to $\sigma_U(c, \varepsilon)$ and
   $A_{c,\varepsilon}(s)\Uparrow$ over $\sigma_U(c, \varepsilon)$ to $\infty$;
(iii) $A_{c,\varepsilon}(s)$ is concave for $s < \sigma_M$ and convex for $s > \sigma_M$.

Note that the critical points of $A_{c,\varepsilon}(s)$ correspond to the two solutions to $\phi(s) = 1/[(1-\varepsilon)c]$. The value of $\sigma_M$, though, does not depend on $c$ or $\varepsilon$. Graphs of a typical function $A_{c,\varepsilon}(\sigma)$ for two different values of $\varepsilon$ are given in Figure 5.

Some further properties, easily verified, are the following:

(a) $\gamma_{b,\varepsilon}$, $\sigma_{b,\varepsilon}$, $\sigma_L(c,\varepsilon)$, $\sigma_U(c,\varepsilon)$ and $A_{c,\varepsilon}(s)$ are continuous in $\varepsilon$;
(b) As $\varepsilon \Uparrow$: $\gamma_{b,\varepsilon} \Downarrow$, $\sigma_{b,\varepsilon} \Uparrow$, $\sigma_L(c,\varepsilon) \Uparrow$, $\sigma_U(c,\varepsilon) \Downarrow$ and $A_{c,\varepsilon}(s) \Downarrow$;
(c) $\gamma_{b,\varepsilon} \leq \sigma_{b,\varepsilon}$ with $\gamma_{b,0} = \sigma_{b,0}$;
(d) if $\gamma < \sigma$, then $A_{c,\varepsilon}(\gamma) - A_{c,\varepsilon}(\sigma)$ is decreasing in $\varepsilon$.

Now, for $1/K < c \leq c_o$,

(A.43) $$\text{if } \sigma_{b,\varepsilon} \leq \sigma_U(c,\varepsilon), \quad \text{then } B_{CM}(\varepsilon) \leq B_S(\varepsilon),$$

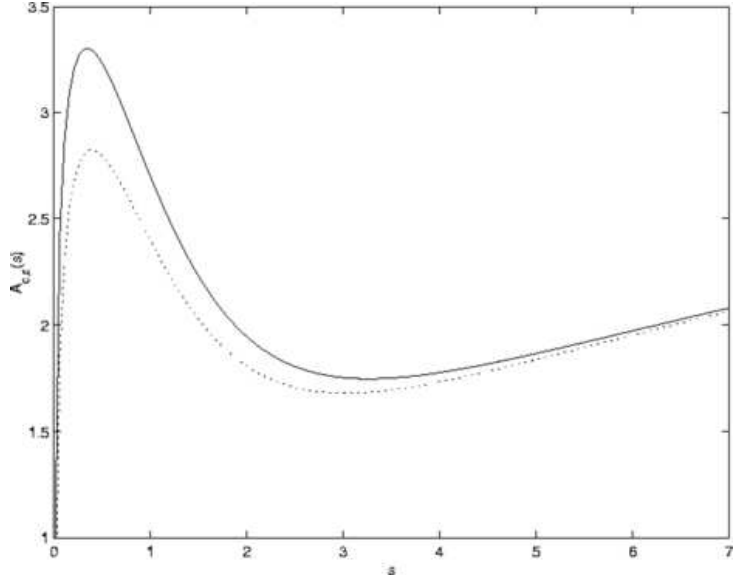

Fig. 5. *Graph of $A_{c,\varepsilon}(\sigma)$.*

since, in this case, $d_c(\varepsilon) \leq 0$. So, to prove Theorem 6.2, it only needs to be shown that

(A.44) $\qquad$ if $\sigma_{b,\varepsilon} > \sigma_U(c,\varepsilon)$, then $A_{c,\varepsilon}(\gamma_{b,\varepsilon}) \leq A_{c,\varepsilon}(\sigma_U(c,\varepsilon))$,

since this implies $d_c(\varepsilon) = 0$ and hence $B_{CM}(\varepsilon) = B_S(\varepsilon)$.

To show (A.44), first note that $\sigma_{b,0} \leq \sigma_M$ since $g(\sigma_{b,0}) = b \geq g(\sigma_M)$. Thus, since $\sigma_{b,\varepsilon}\Uparrow$ and $\sigma_U(c,\varepsilon)\Downarrow$ as $\varepsilon$ increases and both are continuous, there exists an $\varepsilon_b$ such that $\sigma_{b,\varepsilon_b} = \sigma_U(c,\varepsilon_b)$. For any $\varepsilon \leq \varepsilon_b$, it then follows that $\sigma_{b,\varepsilon} \leq \sigma_{b,\varepsilon_b} = \sigma_U(c,\varepsilon_b)$, so to show (A.44), it is only necessary to consider $\varepsilon > \varepsilon_b$.

For $\varepsilon > \varepsilon_b$, we have

$$A_{c,\varepsilon_b}(\gamma_{b,\varepsilon}) \leq A_{c,\varepsilon_b}(\gamma_{b,\varepsilon_b}) \leq A_{c,\varepsilon_b}(\sigma_{b,\varepsilon_b}) = A_{c,\varepsilon_b}(\sigma_U(c,\varepsilon_b)) \leq A_{c,\varepsilon_b}(\sigma_U(c,\varepsilon)).$$

The first inequality follows since $\gamma_{b,\varepsilon_b} \leq \sigma_L(c,\varepsilon_b)$, the second inequality follows from (A.41) and the third inequality follows from (b) since $\sigma_U(c,\varepsilon_b) > \sigma_U(c,\varepsilon) > \sigma_M$. Statement (A.44) then follows from (d) above. $\square$

PROOF OF REMARK 6.1. The remark has already been established for $c > c(\varepsilon)$ and for $c \leq 1/K$. If $c > 1/K$ and $g(\sigma_M) > b$, then $\gamma_{b,0} = \sigma_{b,0} > \sigma_M$. Now, if $\sigma_{b,0} \geq \sigma_U(c,0)$, then since $\sigma_{b,\varepsilon}\Uparrow$ and $\sigma_U(c,0)\Downarrow$ as $\varepsilon$ increases, it follows that $\sigma_{b,\varepsilon} \geq \sigma_U(c,\varepsilon)$ for all $\varepsilon$. This then implies that $d_c(\varepsilon) \geq 0$ and hence that $B_{CM}(\varepsilon) \geq B_S(\varepsilon)$.

On the other hand, if $\sigma_M < \sigma_{b,0} < \sigma_U(c,\varepsilon)$, then by continuity, for sufficiently small $\varepsilon$, we have $\sigma_M < \gamma_{b,\varepsilon} < \sigma_{b,\varepsilon} < \sigma_U(c,\varepsilon)$. This implies that $A_{c,\varepsilon}(\gamma_{b,\varepsilon}) < A_{c,\varepsilon}(\sigma_{b,\varepsilon})$, so by (A.41), $c > c(\varepsilon)$. $\square$

PROOF OF THEOREM 6.3. Note that under the conditions of Theorem 6.2, $B_S(\varepsilon) > B_{CM}(\varepsilon)$ if and only if

$$\sigma_{b,\varepsilon} < \sigma_U(c,\varepsilon) \quad \text{and} \quad A_{c,\varepsilon}(\sigma_{b,\varepsilon}) > A_{c,\varepsilon}(\sigma_U(c,\varepsilon)). \tag{A.45}$$

So, to prove that an $S$-functional is inadmissible, one only needs to establish (A.45) for some $\varepsilon$. First, we will show that the condition $c_o = \lim_{\varepsilon \to 0^+} c(\varepsilon)$ implies that there exists some $\varepsilon$ such that (A.45) holds. Then we will show that (6.29) is enough to guarantee that $c_o = \lim_{\varepsilon \to 0^+} c(\varepsilon)$. Note that by using l'Hôpital's rule, one obtains

$$c(0) = \lim_{\varepsilon \to 0^+} c(\varepsilon) = \frac{1}{\phi(\sigma_{b,0})}. \tag{A.46}$$

Also, note that

$$c > c_1 = \frac{\log(\sigma_K/\sigma_{b,0})}{b - g(\sigma_M)} \iff A_{c,0}(\sigma_{b,0}) > A_{c,0}(\sigma_M). \tag{A.47}$$

Since $\sigma_{b,0} < \sigma_M$, this implies that $c_1 \geq 1/K$ since, otherwise, $A_{c,0}(s)$ would be monotone in $s$. Now, for any $c > c_1$, we then have $\sigma_{b,0} < \sigma_M < \sigma_U(c,0)$ and $A_{c,0}(\sigma_{b,0}) > A_{c,0}(\sigma_M) \geq A_{c,0}(\sigma_U(c,0))$. By continuity, statement (A.45) then follows for sufficiently small $\varepsilon$. Now, we show that $c_1 \leq c(0)$. To show this, note that when $c = c(0)$, we have $\sigma_{b,0} = \sigma_L(c,0)$ and so $A_{c,0}(\sigma_{b,0}) > A_{c,0}(\sigma_M)$. The first part of the proof then follows from (A.47).

Note that the lower bound $c_1$ can be tightened by working with (A.45) directly. In general, it is difficult to use (A.45) to obtain a closed form expression, but it can be used for specific examples.

From (A.46), in the second part of the proof, we need to show that (6.29) implies

$$\varepsilon \, c(\varepsilon) \geq \varepsilon/\phi(\sigma_{b,0}). \tag{A.48}$$

Since equality holds in (A.48) when $\varepsilon = 0$, to show (A.48), it is sufficient to prove that the derivative of the left-hand side is never less than the derivative of the right-hand side, that is [see Equations (A.42) and (A.46)],

$$\frac{1}{(1-\varepsilon)^2}\left\{\frac{1-b}{\phi(\sigma_{b,\varepsilon})} + \frac{b}{\phi(\gamma_{b,\varepsilon})}\right\} \geq \frac{1}{\phi(\sigma_{b,0})}. \tag{A.49}$$

Recall that we are assuming $g(\sigma_M) < b = g(\sigma_{b,0})$ or, equivalently, that $\sigma_{b,0} < \sigma_M$. This implies $\phi(\gamma_{b,\varepsilon}) < \phi(\sigma_{b,0})$ and after some simple algebraic manipulations, we note that (A.49) holds if

$$a_{b,\varepsilon}\phi(\sigma_{b,\varepsilon}) \leq \phi(\sigma_{b,0}), \tag{A.50}$$

where $a_{b,\varepsilon} = [(1-\varepsilon)^2 - b]/(1-b)$.

Since $\sigma_{b,\varepsilon}$ is increasing in $\varepsilon$, it follows that $\phi(\sigma_{b,\varepsilon})$ is decreasing in $\varepsilon$ whenever $\sigma_{b,\varepsilon} \geq \sigma_M$ and hence that if (A.50) holds for $\sigma_{b,\varepsilon} = \sigma_M$, then it also holds for $\sigma_{b,\varepsilon} \geq \sigma_M$. Thus, it is sufficient to show that (A.50) holds for $\sigma_{b,\varepsilon} \leq \sigma_M$ or, equivalently, for

$$\varepsilon \leq \varepsilon_M \doteq \frac{b - g(\sigma_M)}{1 - g(\sigma_M)}.$$

Given that $g(s)$ is convex, we have $-g'(\sigma_{b,\varepsilon}) \leq -g'(\sigma_{b,0})$, so (A.50) holds if $a_{b,\varepsilon}\ \sigma_{b,\varepsilon} \leq \sigma_{b,0}$. Since $g(s)$ is also nonincreasing, this is equivalent to

$$g(a_{b,\varepsilon}\ \sigma_{b,\varepsilon}) \geq b. \tag{A.51}$$

Thus, the theorem is proved if (A.51) holds for $\varepsilon \leq \varepsilon_K$. By the convexity of $g(s)$, for $\varepsilon \leq \varepsilon_K$, we have

$$g(a_{b,\varepsilon}\sigma_{b,\varepsilon}) \geq g(\sigma_{b,\varepsilon}) + (a_{b,\varepsilon} - 1)\sigma_{b,\varepsilon}g'(\sigma_{b,\varepsilon})$$
$$= \frac{b-\varepsilon}{1-\varepsilon} + \frac{\varepsilon(2-\varepsilon)}{1-b}\phi(\sigma_{b,\varepsilon}) \geq \frac{b-\varepsilon}{1-\varepsilon} + \frac{\varepsilon(2-\varepsilon)}{1-b}\phi(\sigma_{b,0}).$$

The last term is $\geq b$ if and only if

$$\phi(\sigma_{b,0}) \geq \frac{(1-b)^2}{(1-\varepsilon)(2-\varepsilon)}. \tag{A.52}$$

Note that if (A.52) holds for $\varepsilon = \varepsilon_M$, then it holds for all $\varepsilon \leq \varepsilon_M$. With $\varepsilon = \varepsilon_M$, though, (A.52) corresponds to the bound (6.29). This completes the proof. $\square$

J. R. BERRENDERO
DEPARTAMENTO DE MATEMÁTICAS
UNIVERSIDAD AUTÓNOMA DE MADRID
28049 MADRID
SPAIN
E-MAIL: joser.berrendero@uam.es

B. V. M. MENDES
INSTITUTO DE MATEMÁTICA / COPPEAD
UNIVERSIDADE FEDERAL DO RIO DE JANEIRO
21941-918 RJ
BRAZIL
E-MAIL: beatriz@im.ufrj.br

D. TYLER
DEPARTMENT OF STATISTICS
HILL CENTER, BUSCH CAMPUS
RUTGERS UNIVERSITY
PISCATAWAY, NEW JERSEY 08854
USA
E-MAIL: dtyler@rci.rutgers.edu